\documentclass{article}
\usepackage{amssymb,amsmath,amsthm}
\usepackage{latexsym}

\usepackage[dvips]{epsfig,graphics}
\input{epsf}

\newtheorem{theorem}{Theorem}

\newtheorem{definition}[theorem]{Definition}

\newtheorem{lemma}[theorem]{Lemma}

\newtheorem{proposition}[theorem]{Proposition}
\newtheorem{remark}[theorem]{Remark}

\input{tcilatex}

\begin{document}

\title{The coefficients of the HOMFLYPT and the Kauffman polynomials are pointwise limits
of Vassiliev invariants.}
\author{Laure Helme-Guizon \\
%EndAName
George Washington University}
\maketitle

\begin{abstract}
The Vassiliev conjecture states that the Vassiliev invariants are dense in the space of all numerical link
invariants in the sense that any link invariant is a pointwise limit of Vassiliev invariants.

In this article, we prove that the Vassiliev conjecture holds in the case of  the coefficients of the HOMFLY and
the Kauffman polynomials.
\end{abstract}

\tableofcontents

\bigskip

A well-known conjecture in the theory of Vassiliev invariants, called
\textit{Vassiliev conjecture}, is that these invariants are dense in the
space of all numerical link invariants. This was posed as a problem in
[BL93] as follows: Given any numerical link invariant $f:\mathcal{L}%
\rightarrow \mathbb{Q},$ does there exist a sequence of Vassiliev invariants
$\{v_{n}^{f}:\mathcal{L}\rightarrow \mathbb{Q}\}_{n\in \mathbb{N}}$ such
that, for any fixed $L$, $\underset{n\rightarrow \infty }{\lim }%
v_{n}^{f}(L)=f(L)$ ?

In other words, is any link invariant a pointwise limit of Vassiliev
invariants?

\medskip

Vassiliev conjecture arouses interest because, if it holds, it would imply
that the Vassiliev invariants can distinguish any two links, as can be seen
using the approximation by Vassiliev invariants of the link invariant $%
\delta _{L_{1}}$ defined by $\delta _{L_{1}}(L)$ $=1$ if $L=L_{1},$ $\delta
_{L_{1}}(L)=0$ otherwise.\medskip

It has been proven in [KR00] that the coefficients of the Jones polynomial
of a knot are pointwise limit of Vassiliev invariants. As noted in the
article, the result carries over with slight modifications to links.

But when one tries to extend these results to the coefficients of the
HOMFLYPT and the Kauffman 2-variable polynomial, the method used to achieve
the previous result fails. Indeed, the idea was to produce some Vassiliev
invariants that are linear combinations of the coefficients of the Jones
polynomial, through a change of variable in the polynomial. Then one would
solve back using a matrix.

In the case of 2-variable polynomials, one has a double sum so one cannot
solve back the same way and the standard change of variable gives a pretty
complicated formula for the Vassiliev invariants. The approach of this paper
differs from the previous one in three ways : First, we will use a non
standard change of variable to produce Vassiliev invariants, which leads to
easier formulas. This is explained in subsection $2.2$. Second, we will
introduce intermediate variables in order to get sums on one index only.
These intermediate variables may not be Vassiliev invariants, but their
restrictions to links with a fixed number of components are, as explained in
subsection $2.3$. So, and this is the third difference, we need to prove an
equivalent but seemingly weaker formulation of the Vassiliev conjecture
which says that working with links with a fixed number of components is
enough. This is explained in section $1$.

All these steps together reduce the problem to a point where the method \ of
[KR00] applies. We will supply an alternative proof for this last step,
which uses complex analysis. These two proofs can be found in subsection $2.4
$.

In section $3$, we explain how a very similar argument can be used in the
case of the Kauffman 2-variable polynomial.\medskip

All this will allow us to show that Vassiliev conjecture also holds for the
coefficients of the HOMFLYPT and the Kauffman 2-variable polynomial.\medskip

\section{Reducing the problem: the weak Vassiliev conjecture}

In order to prove this result, we will use a reformulation of\ Vassiliev
conjecture, which we\ call the \textit{weak Vassiliev conjecture}, that
seems to be weaker but which we will prove to be equivalent to Vassiliev
conjecture.\medskip

Let $\mathcal{L}$ be the set of all oriented links and let $\mathcal{L}%
^{(\mu )}$ be the set of oriented links with $\mu $ components. We can define what it means for a link invariant
defined on $\mathcal{L}^{(\mu )}$ only to be a Vassiliev invariant the following way: A link with
self-intersections is said to be a $\mu $ \textit{component link} if by replacing all the singular crossings by
positive or negative crossing we get a $\mu $ component link. This definition makes sense because switching
crossings does not affect the number of components. We denote $X_{i}^{(\mu )} $ the set of $\mu $ component
links with $i$ self-intersections. Let $ v^{\mu }:\mathcal{L}^{(\mu )}\rightarrow \mathbb{Q}$ be a link
invariant defined on the set of oriented links with $\mu $ components. We can extend it by induction to all the
$X_{i}^{(\mu )}$'s using \textit{Vassiliev relation} $v^{\mu }(L_{\times })=v^{\mu }(L_{+})-v^{\mu }(L_{-})$
where\ the diagrams of $L_{\times },L_{+}$ and $L_{-}$ are the same except in a neighborhood of a self
intersection, denoted by $\times ,$ which is replaced by a positive crossing in $L_{+}$ and by a negative
crossing in $L_{-}$. If $ v^{\mu }$ vanishes on all link with $n+1$ singular crossings, $v^{\mu }$ is said to be
a \textit{Vassiliev invariant of order at most }$\mathit{n}$ \textit{.}\medskip

We are now ready to state the \textit{weak Vassiliev conjecture} : For any link invariant
$f:\mathcal{L}\rightarrow \mathbb{Q}$\textbf{, }for any $\mu \in \mathbb{N}^{\ast }(=\mathbb{N\setminus }\left\{
0\right\} )$, there exists a sequence of Vassiliev invariants $\{v_{n}^{f,\mu }:\mathcal{L} ^{(\mu )}\rightarrow
\mathbb{Q}\}_{n\in \mathbb{N}}$ such that, for any fixed $L\in \mathcal{L}^{(\mu )},\underset{n\rightarrow
\infty }{\lim } v_{n}^{f,\mu }(L)=f(L).$

\begin{proposition}
\label{WVCe}If all the restrictions of a link invariant to links with a
fixed number of component are pointwise limits of Vassiliev invariants then
the link invariant itself is a pointwise limit of Vassiliev invariants. In
other words, the weak Vassiliev conjecture and\ Vassiliev conjecture are
equivalent.
\end{proposition}

\begin{proof}
Let $\ f:\mathcal{L}\rightarrow \mathbb{Q}$ be a\textbf{\ }link invariant such that\textbf{\ }for any $\mu \in
N^{\ast }$, there exists a sequence of Vassiliev invariants $\{v_{n}^{\mu }:\mathcal{L}^{(\mu )}\rightarrow
\mathbb{ Q}\}_{n\in \mathbb{N}}$ satisfying the condition that, for any fixed $L\in \mathcal{L}^{(\mu
)},\underset{n\rightarrow \infty }{\lim }v_{n}^{f,\mu }(L)=f(L).$ Let $\{w_{n}^{\mu }:\mathcal{L}^{(\mu
)}\rightarrow \mathbb{Q} \}_{n\in \mathbb{N}}$ be defined by $w_{n}^{\mu }=v_{n}^{\mu }$ if $n\geq \mu $ and $0$
otherwise. Each $w_{n}^{\mu }$ is still a Vassiliev invariant. For any given value of $\mu ,$
$\underset{n\rightarrow \infty }{\lim } w_{n}^{\mu }(L)=f(L)$ because perturbing a finite number of terms in a
sequence doesn't change the limit.\medskip

Let $w_{n}:\mathcal{L}\rightarrow \mathbb{Q}$ be the\textbf{\ }link invariant defined by $w_{n}(L)=w_{n}^{\mu
}(L)$ if $L\in \mathcal{L}^{(\mu )}.$ $w_{n}$ is a Vassiliev invariant with order$\ =$ $\underset{\mu \in
N^{\ast }}{\sup }\left\{ order\left( w_{n}^{\mu }\right) \right\} =$ $ \underset{\mu \leq n}{\max }\left\{
order\left( w_{n}^{\mu }\right) \right\} $ which is well defined because\ the second set is finite.

And obviously, $\underset{n\rightarrow \infty }{\lim }w_{n}(L)=f(L).$
\end{proof}

\section{Approximation of the coefficients of the HOMFLYPT polynomial by
Vassiliev invariants}

\subsection{Notations}

Let $P_{L}(v,z)$ denote the HOMFLYPT polynomial of a link $L$. The version
of the HOMFLYPT polynomial we are working with is the one defined by :

\begin{enumerate}
\item  $v^{-1}P_{L_{+}}(v,z)-vP_{L_{-}}(v,z)$ $=zP_{L_{0}}(v,z).$

\item  $P_{T_{1}}(v,z)$ $=1$ where $T_{1}$ is the unknot of one
component.\medskip
\end{enumerate}

Suppose the HOMFLYPT polynomial of a link $L$ is

$P_{L}(v,k)=\underset{-k_{1}\leqslant k\leqslant k_{2}}{\sum }\ \underset{ -j_{1}\leqslant j\leqslant
j_{2}}{\sum }a_{kj}(L)v^{k}z^{j}$ where $k_{1},$ $ k_{2},j_{1},$ $j_{2}\in \mathbb{N}$

We call the \textit{degree} of such a Laurent polynomial $ d=max\{k_{1},k_{2},j_{1},j_{2}\}$.\medskip

\begin{lemma}
Let $\mu $ be the number of components of a link $L.$ Its HOMFLYPT polynomial is a polynomial (not a Laurent
polynomial) in the variables $ v,v^{-1}\ ,$ $z$ and $\frac{v-v^{-1}}{z}$, the maximum possible power of $
\frac{v-v^{-1}}{z}$ being $\mu -1.$ \label{max pow 1/z}
\end{lemma}

\begin{proof}
The unknotting number of a link diagram $D$ is the minimum number of
crossings one needs to switch to get the trivial link. Note that it is not
necessarily equal to unknotting number of the link, which takes into account
all the diagrams.

The lemma can be proven by induction on $(c(D),u(D))$ where $D$ is a diagram
of the link $L$, $c(D)$ is the number of crossing of $D,$ $u(D)$ is the
unknotting number of $D$ and the order on $\mathbb{N}\times \mathbb{N}$ is
the lexicographic order.\medskip
\end{proof}

\subsection{Step 1 : Get Vassiliev invariant through a power series expansion
}

It is proved in [BL93] that the change of variables $v=e^{-\frac{N}{2}x}$ \
and $z=e^{\frac{1}{2}x}-e^{-\frac{1}{2}x}$ in $P_{L}(v,z)$ gives a power
series expression for $P_{L}(v,z)$ in variable $x$ whose coefficients are
Vassiliev invariants.

This change of variables has become the standard one and is the one that was
used in [KR00]. But in the case of the HOMFLYPT polynomial, it leads to very
complicated formulas when it comes to expressing the Vassiliev invariants in
terms of the coefficients of the polynomial.

It happens that we get much easier computations using the change of
variables explained in the following lemma.

\begin{lemma}
Let $W_{N,x}(L)$ be the polynomial obtained from the HOMFLYPT polynomial of $ L$ by setting $v:=e^{Nx}$ and
$z:=x.$ After power series expansion, each $ w_{Nq}(L)$ in $W_{N,x}(L)=\underset{q=0}{\overset{\infty }{\sum }}
w_{Nq}(L)x^{q}$ is a Vassiliev invariant of order (at most) $q$, for all $N$ $\in \mathbb{Z}$.\label{Chg Var HOM
VI}
\end{lemma}

Notice that lemma (\ref{max pow 1/z}) ensures that they are no negative powers of x after power series expansion
because each of the variables $ v,v^{-1}\ ,$ $z$ and $\frac{v-v^{-1}}{z}$ has a power series expansion.\medskip

\begin{proof}
The proof is similar to the one used in [BL93] to show that if we let $v=e^{- \frac{N}{2}x}$ \ and
$z=e^{\frac{1}{2}x}-e^{-\frac{1}{2}x}$ in $P_{L}(v,z)$\ we get, after power series expansion,
$\widetilde{W}_{K}(N,x)=\underset{q=0}{ \overset{\infty }{\sum }}\widetilde{w}_{Nq}(K)x^{q}$\ \ where
$\widetilde{w}
_{Nq}(K)$ is a Vassiliev invariant of order (at most) $q$, for all $N$%
.\medskip
\end{proof}

By lemma (\ref{max pow 1/z}), the HOMFLYPT polynomial of a link $L$ with $%
\mu $ components has no power of $z$ with exponent less than $-\mu +1,$ hence we can rewrite it \
$P_{L}(v,k)=\underset{-d\leqslant k\leqslant d}{ \sum }$ $\underset{-\mu +1\leqslant j\leqslant d}{\sum
}a_{kj}(L)v^{k}z^{j} \label{Hkmu}$ where $d$ is the degree of the HOMFLYPT polynomial of $L$, with the
convention $a_{k,j}=o$ if $\left| k\right| >d$ or $\left| j\right|
>d.$\medskip

Let $W_{N,x}(L)$ be the polynomial obtained from the HOMFLYPT polynomial of $ L$ by setting $v:=e^{Nx}$ and
$z:=x.$ Let us compute the power series expansion :

$W_{N,x}(L)=\underset{-d\leqslant k\leqslant d}{\sum }\ \underset{-\mu +1\leqslant j\leqslant d}{\sum
}a_{kj}(L)e^{Nkx}x^{j}=\underset{-d\leqslant k\leqslant d}{\sum }\ \underset{-\mu +1\leqslant j\leqslant d}{\sum
} a_{kj}(L)\ \underset{s=0}{\overset{\infty }{\sum }}\frac{N^{s}k^{s}x^{s}}{s!} x^{j}.$

Let $q:=s+j,$ we get $W_{N,x}(L)=$\ $\underset{-\mu +1\leqslant j\leqslant d }{\sum }\ \underset{-d\leqslant
k\leqslant d}{\sum }\ a_{kj}(L)\underset{q=j }{\overset{\infty }{\sum }}\frac{N^{q-j}k^{q-j}}{(q-j)!}x^{q}.$

Exchanging the sums over $k$ and $q$ yields $W_{N,x}(L)=$\ $\underset{-\mu +1\leqslant j\leqslant d}{\sum }\
\underset{q=j}{\overset{\infty }{\sum }}\ \underset{-d\leqslant k\leqslant d}{\sum
}a_{kj}(L)\frac{N^{q-j}k^{q-j}}{ (q-j)!}x^{q}.$

Exchanging (carefully) the sums over $j$ and $q$ yields

$W_{N,x}(L)=$\ $\underset{q=-\mu +1}{\overset{\infty }{\sum }}\underset{ w_{N,q}(L)}{\ \underbrace{\left(
\underset{-\mu +1\leqslant j\leqslant \min (q,d)}{\sum }\ \underset{-d\leqslant k\leqslant d}{\sum
}a_{kj}(L)\frac{ N^{q-j}k^{q-j}}{(q-j)!}\right) }}x^{q}$

By Lemma $\left( \ref{Chg Var HOM VI}\right) $, we know that the $w_{Nq}(L)$
are Vassiliev invariants. The above formula shows that they are related to
the initial coefficients of the HOMFLYPT polynomial $a_{kj}(L)$ by

\begin{equation*}
w_{Nq}(L)=\sum_{-\mu +1\leqslant j\leqslant \min (q,d)}N^{q-j}\sum_{-d\leqslant k\leqslant
d}a_{k,j}(L)\frac{k^{q-j}}{(q-j)!} \text{ for all }q,N\in \mathbb{Z}\text{ with }q\geq -\mu +1
\end{equation*}

After the change of variable $p=q-j,$ we get :

\begin{equation}
w_{Nq}(L)=\sum_{0\leqslant p\leqslant q+\mu -1}N^{p}\sum_{-d\leqslant k\leqslant
d}a_{k,q-p}(L)\frac{k^{p}}{p!}\text{ for all }q,N\in \mathbb{Z} \text{ with }q\geq -\mu +1  \label{wNq}
\end{equation}

Our goal is now to show that the coefficients $a_{kj}$ of the HOMFLYPT
polynomial satisfy the weak Vassiliev conjecture.\medskip

Note that because of the double sum, the method used in [KR00] to prove that
the coefficients of the Jones polynomial of a link are limits of Vassiliev
invariants cannot be applied to prove that the coefficients of the HOMFLYPT
polynomial of a link are limits of Vassiliev invariants. Therefore, we will
introduce new variables:\medskip

\subsection{Step 2: The Intermediate variables $B_{mj}^{\protect\mu }$ are
Vassiliev invariants.}

\begin{definition}
Let $B_{mj}(L)=\underset{k=-d}{\overset{d}{\sum }}a_{kj}(L)\frac{k^{m}}{m!}$
for all $m$ $\in \mathbb{N},$ $j$ $\in \mathbb{Z}$.
\end{definition}

\begin{definition}
Let $B_{mj}^{\mu }$ (resp. $a_{kj}^{\mu })$ be the restriction of $B_{mj}$
(resp. $a_{kj})$ to the links of $\mu $ components.
\end{definition}

We have $B_{mj}^{\mu }$ $=\underset{k=-d}{\overset{d}{\sum }}a_{kj}^{\mu }(L) \frac{k^{m}}{m!}.$

\begin{proposition}
The $B_{mj}^{\mu }$ are Vassiliev invariants, for all $\mu \in \mathbb{N} ^{\ast },$ $m$ $\in \mathbb{N},$ $j$
$\in \mathbb{Z}\label{BmjVI}$
\end{proposition}

\begin{proof}
Let $\mu \in \mathbb{N}^{\ast },m\in \mathbb{N},$ $j$ $\in \mathbb{Z}$.

$\blacktriangle $ \ \emph{Case 1 : Assume }$j<-\mu +1.$

We already noticed that $P_{L}(v,k)$ can be written $P_{L}(v,k)=\underset{ -d\leqslant k\leqslant d}{\sum }$
$\underset{-\mu +1\leqslant j\leqslant d}{ \sum }a_{kj}(L)v^{k}z^{j},$ so $a_{kj}(L)=0$ \ for all $k$ when$\
j<-\mu +1, \emph{\ }$hence $B_{mj}^{\mu }(L)=0$ for all $L\in \mathcal{L}^{(\mu )}.$ Therefore, it is a
Vassiliev invariant.\medskip

$\blacktriangle $ \ \emph{Case 2 : Assume }$j\geq -\mu +1.$

Let $q:=m+j.$ $m\in \mathbb{N}$ so $q=m+j\geq -\mu +1.$ By formula (\ref{wNq} ), the $w_{Nq}^{\mu }$ can be
expressed in terms of the $B_{mj}^{\mu }$'s by
$\ \ \ $%
\begin{equation}
w_{Nq}^{\mu }(L)=\sum_{0\leqslant p\leqslant q+\mu -1}N^{p}\text{ } B_{p,q-p}^{\mu }(L)\quad \text{for all N}\in
\mathbb{Z}\text{, for all }L\in \mathcal{L}^{(\mu )}  \label{wNq in terms of Bmj}
\end{equation}

By letting $N=1,2,...,q+\mu $ , we get a set of equalities that can be
summarized in a matrix equality,

$\left(
\begin{array}{cccc}
1 & 1 & ... & 1 \\
1 & (2)^{1} & ... & (2)^{^{q+\mu -1}} \\
1 & (3)^{1} & ... & (3)^{^{q+\mu -1}} \\
\vdots & \vdots &  & \vdots \\
1 & (q+\mu )^{1} & ... & (q+\mu )^{^{q+\mu -1}}
\end{array}
\right) \cdot \left(
\begin{array}{c}
B_{0,q}^{\mu }(L) \\
B_{1,q-1}^{\mu }(L) \\
\vdots \\
B_{q+\mu -1,-\mu +1}^{\mu }(L)
\end{array}
\right) =\left(
\begin{array}{c}
w_{1,q}^{\mu }(L) \\
w_{2,q}^{\mu }(L) \\
w_{3,q}^{\mu }(L) \\
\vdots \\
w_{q+\mu ,q}^{\mu }(L)
\end{array}
\right) $

Let $n$ $\geq $ $q+\mu .$ Let $A_{n}=\left(
\begin{array}{cccc}
1 & 1 & ... & 1 \\
1 & (2)^{1} & ... & (2)^{n-1} \\
1 & (3)^{1} & ... & (3)^{n-1} \\
\vdots & \vdots &  & \vdots \\
1 & (n)^{1} & ... & (n)^{n-1}
\end{array}
\right) $. \ As seen in case 1, $B_{m\widetilde{j}}^{\mu }=0$ whenever $ \widetilde{j}\leq -\mu ,$ so$\left(
\begin{array}{c}
B_{0,q}^{\mu }(L) \\
B_{1,q-1}^{\mu }(L) \\
\vdots \\
B_{n-1,q-n+1}^{\mu }(L)
\end{array}
\right) =\left(
\begin{array}{c}
B_{0,q}^{\mu }(L) \\
B_{1,q-1}^{\mu }(L) \\
\vdots \\
B_{q+\mu -1,-\mu +1}^{\mu }(L) \\
0 \\
\vdots \\
0
\end{array}
\right) $

Hence, for all $n$ $\geq $ $q+\mu ,$ $\ A_{n}\cdot \left(
\begin{array}{c}
B_{0,q}^{\mu }(L) \\
B_{1,q-1}^{\mu }(L) \\
\vdots \\
B_{n-1,q-n+1}^{\mu }(L)
\end{array}
\right) =\left(
\begin{array}{c}
w_{1,q}^{\mu }(L) \\
w_{2,q}^{\mu }(L) \\
w_{3,q}^{\mu }(L) \\
\vdots \\
w_{n,q}^{\mu }(L)
\end{array}
\right) $

$A_{n}$ is a $n\times n$ Vandermonde matrix with distinct parameters, thus $A_{n}$ is invertible.

Hence, for all $n$ $\geq $ $q+\mu ,\left(
\begin{array}{c}
B_{0,q}^{\mu }(L) \\
B_{1,q-1}^{\mu }(L) \\
\vdots \\
B_{n-1,q-n+1}^{\mu }(L)
\end{array}
\right) =\left( A_{n}\right) ^{-1}\cdot \left(
\begin{array}{c}
w_{1,q}^{\mu }(L) \\
w_{2,q}^{\mu }(L) \\
w_{3,q}^{\mu }(L) \\
\vdots \\
w_{n,q}^{\mu }(L)
\end{array}
\right) $

Since the coefficients and the size of $\left( A_{n}\right) ^{-1}$do not depend on $L$ (remember that $\mu $ is
fixed), each $B_{i,q-i}^{\mu },$ $ 0\leq i\leq n-1$ is a linear combination of the Vassiliev invariants $
w_{i,q}^{\mu }$'s, so each is a Vassiliev invariant. Letting $n:=\max (m+1,$ $q+\mu ),$ we get that $B_{mj}^{\mu
}$ is a Vassiliev invariant.
\end{proof}

\subsection{Last step : The initial coefficient are limits of Vassiliev
invariants.}

We are now ready to prove our main result, namely that the $a_{kj}(L)$'s,
the initial coefficient of the HOMFLYPT polynomial are limits of Vassiliev
invariants. We know that it suffices to show that they satisfy the weak
Vassiliev conjecture.

At this point, given that the $B_{mj}^{\mu }$'s are Vassiliev invariant, the
method\ exposed by Y. Rong and I. Kofman in [KR00] can be implemented. We
will therefore give a first proof based on this method.

The idea that this result could probably also be proven using complex
analysis techniques was suggested by Jo$\widetilde{a}$o Faria Martins and
lead to a second proof, exposed thereafter.

\begin{proposition}
The coefficients $a_{kj}(L)$ of the HOMFLYPT polynomial satisfy the weak
Vassiliev conjecture i.e., their restrictions to links with a fixed number
of components is a pointwise limit of Vassiliev invariants.
\end{proposition}

\begin{proof}
1\textit{, Using Linear Algebra}

$\blacktriangle $ Let $j$ $\in \mathbb{Z}$. By letting $m=0,1,2,...$ \ in $
B_{mj}(L)=\underset{k=-d}{\overset{d}{\sum }}a_{kj}(L)\frac{k^{m}}{m!}$, we get a system of equations that can
be expressed with an infinite matrix$.$.

\begin{equation}
\left(
\begin{array}{ccccccc}
... & 1 & 1 & 1 & 1 & 1 & ... \\
... & -2 & -1 & 0 & 1 & 2 & ... \\
... & \left( -2\right) ^{2} & \left( -1\right) ^{2} & \left( 0\right) ^{2} &
\left( 1\right) ^{2} & \left( 2\right) ^{2} & ... \\
... & \left( -2\right) ^{3} & \left( -1\right) ^{3} & \left( 0\right) ^{3} &
\left( 1\right) ^{3} & \left( 2\right) ^{3} & ... \\
... & \left( -2\right) ^{4} & \left( -1\right) ^{4} & \left( 0\right) ^{4} &
\left( 1\right) ^{4} & \left( 2\right) ^{4} & ... \\
& \vdots & \vdots & \vdots & \vdots & \vdots &
\end{array}
\right) \cdot \left(
\begin{array}{c}
\vdots \\
a_{-2,j}(L) \\
a_{-1,j}(L) \\
a_{0,j}(L) \\
a_{1,j}(L) \\
a_{2,j}(L) \\
\vdots
\end{array}
\right) =\left(
\begin{array}{c}
0!B_{0j}(L) \\
1!B_{1j}(L) \\
2!B_{2j}(L) \\
3!B_{3j}(L) \\
\vdots
\end{array}
\right)  \label{infMat}
\end{equation}

with the convention $a_{k,j}=o$ if $\left| k\right| >d$ or $\left| j\right|
>d$

We use an infinite matrix because we do not want its size to depend on the
knot and the degree $d$ does.

Let $M_{n}$ be the finite $(2n+1)\times (2n+1)$ matrix extracted from the
above infinite matrix and defined by
\begin{equation*}
M_{n}=\left(
\begin{array}{ccccccc}
1 & ... & 1 & 1 & 1 & ... & 1 \\
-n & ... & -1 & 0 & 1 & ... & n \\
\left( -n\right) ^{2} & ... & \left( -1\right) ^{2} & \left( 0\right) ^{2} &
\left( 1\right) ^{2} & ... & \left( n\right) ^{2} \\
\left( -n\right) ^{3} & ... & \left( -1\right) ^{3} & \left( 0\right) ^{3} &
\left( 1\right) ^{3} & ... & \left( n\right) ^{3} \\
\vdots &  & \vdots & \vdots & \vdots &  & \vdots \\
\left( -n\right) ^{2n} & .. & \left( -1\right) ^{2n} & \left( 0\right) ^{2n}
& \left( 1\right) ^{2n} & .. & \left( n\right) ^{2n}
\end{array}
\right)
\end{equation*}
The matrix $M_{n}$ is a Vandermonde matrix with distinct parameters, and
thus is invertible.

For all $n\in N$, the linear equation $M_{n}\;X_{2n+1}=\left(
\begin{array}{c}
0!B_{0j}(L) \\
1!B_{1j}(L) \\
2!B_{2j}(L) \\
\vdots \\
(2n)!B_{2n,j}(L)
\end{array}
\right) $ has a unique solution $s^{n,j}(L)=\left(
\begin{array}{c}
s_{-n}^{n,j}(L) \\
s_{-n+1}^{n,j}(L) \\
... \\
s_{0}^{n,j}(L) \\
... \\
s_{n}^{n,j}(L)
\end{array}
\right) =\left( M_{n}\right) ^{-1}\cdot \left(
\begin{array}{c}
0!B_{0j}(L) \\
1!B_{1j}(L) \\
2!B_{2j}(L) \\
\vdots \\
(2n)!B_{2n,j}(L)
\end{array}
\right) $

Let $\mu \in N^{\ast }.$ We now restrict all the link invariants to links of
$\mu $ components, and we get :

$s^{n,j,\mu }(L)=\left(
\begin{array}{c}
s_{-n}^{n,j,\mu },(L) \\
s_{-n+1}^{n,j,\mu }(L) \\
... \\
s_{0}^{n,j,\mu }(L) \\
... \\
s_{n}^{n,j,\mu }(L)
\end{array}
\right) =\left( M_{n}\right) ^{-1}\cdot \left(
\begin{array}{c}
0!B_{0j}^{\mu }(L) \\
1!B_{1j}^{\mu }(L) \\
2!B_{2j}^{\mu }(L) \\
\vdots \\
(2n)!B_{2n,j}^{\mu }(L)
\end{array}
\right) $

Since the coefficients and the size of $\left( M_{n}\right) ^{-1}$do not depend on $L$, each $s_{k}^{n,j,\mu
}(L)$ is a linear combination of the Vassiliev invariants $B_{mj}^{\mu }(L)$'s, so each $s_{k}^{n,j,\mu }(L)$ is
itself a Vassiliev invariant for all $k$ such that $-n\leqslant k\leqslant n$ .\medskip

$\blacktriangle $ We are now ready to prove that for any fixed $\mu ,j,k$ and any fixed $L\in \mathcal{L}^{(\mu
)}$, $\underset{n\rightarrow \infty }{ \lim }s_{k}^{n,j,\mu }(L)=a_{kj}^{\mu }(L).$ This follows directly from
the following more general statement (in which the number of components is no longer fixed): For any fixed $j,k$
and any fixed $L\in \mathcal{L}$, $\underset{n\rightarrow \infty }{\lim }s_{k}^{n,j}(L)=a_{kj}(L).$

Let $L$ be a link. Let $k$, $j$ $\in \mathbb{Z}$. Let $n\geq d$ where $d$ is
the degree of HOMFLYPT polynomial of $K.$

We know by (\ref{infMat}) that $M_{n}\left(
\begin{array}{c}
0 \\
\vdots \\
0 \\
a_{-d,j}(L) \\
\vdots \\
a_{0,j}(L) \\
\vdots \\
a_{d,j}(L) \\
0 \\
\vdots \\
0
\end{array}
\right) =\left(
\begin{array}{c}
0!B_{0j}(L) \\
1!B_{1j}(L) \\
2!B_{2j}(L) \\
\vdots \\
(2n)!B_{2n,j}(L)
\end{array}
\right) .$ By unicity of the solution of this matrix equation, we get $ s_{k}^{n,j}(L)=a_{kj}(L)$ for all $n\geq
d.$ Hence the sequence $\left\{ s_{k}^{n,j}(L)\right\} _{n\geq 1}$becomes stationary so $\underset{ n\rightarrow
\infty }{\lim }s_{k}^{n,j}(L)=a_{kj}(L)$.

If we now restrict all the link invariants to links of $\mu $ components, we
get $\underset{n\rightarrow \infty }{\lim }s_{k}^{n,j,\mu }(L)=a_{kj}^{\mu
}(L)$. Since the $s_{k}^{n,j,\mu }$'s are Vassiliev invariants, the
coefficients $a_{kj}^{\mu }$ of the HOMFLYPT polynomial are limits of
Vassiliev invariants as announced.\bigskip
\end{proof}

\begin{proof}
2\textit{: Using Complex analysis}

For any link $L$ and any $j$ $\in \mathbb{Z}$, define $f_{L,j}:$ $\mathbb{C} ^{\ast }\rightarrow \mathbb{C}$ by
$f_{L,j}(z):=\underset{\ \ k=-d}{\overset{ d}{\sum }}a_{kj}(L)z^{k}$ where the $a_{kj}$'s are the coefficients
of the
HOMFLYPT polynomial written \ $P_{L}(v,z)=\underset{-d\leqslant k\leqslant d%
}{\sum }$ $\underset{-\mu +1\leqslant j\leqslant d}{\sum } a_{kj}(L)v^{k}z^{j}.$

$\blacktriangle $ $f_{L,j}(e^{x})=\underset{k=-d}{\overset{d}{\sum }}
a_{kj}(L)e^{kx}=\underset{k=-d}{\overset{d}{\sum }}a_{kj}(L)\underset{m=0}{ \overset{\infty }{\sum
}}\frac{k^{m}x^{m}}{m!}=\underset{m=0}{\overset{ \infty }{\sum }}\left( \underset{k=-d}{\overset{d}{\sum
}}a_{kj}(L)\frac{
k^{m}}{m!}\right) x^{m}$%
\begin{equation}
\text{so }f_{L,j}(e^{x})=\underset{m=0}{\overset{\infty }{\sum }} B_{mj}(L)x^{m}.  \label{f(ex)}
\end{equation}

Let $n$ $\in \mathbb{Z}$. From $\ f_{L,j}(z)=\underset{k=-d}{\overset{d}{ \sum }}a_{kj}(L)z^{k},$ we get, by a
standard result in complex analysis, $ a_{nj}(L)=\frac{1}{2\pi i}\int_{\Gamma ^{+}}\frac{\
f_{L,j}(z)}{z^{n+1}}dz$ where $\Gamma ^{+}$ is the unit circle travelled around once in the positive direction.

$\blacktriangle $ Let's compute this integral:\ $2\pi i\cdot a_{nj}(L)=\int_{\Gamma ^{+}}\frac{\
f_{L,j}(z)}{z^{n+1}}dz=\int_{0}^{2\pi } \frac{\ f_{L,j}(e^{it})}{e^{it(n+1)}}ie^{it}dt$ \ Using (\ref{f(ex)}),$\
2\pi i\cdot a_{nj}(L)=\int_{0}^{2\pi }\frac{\underset{m=0}{\overset{\infty }{ \sum
}}B_{mj}(L)(it)^{m}}{e^{it(n+1)}}ie^{it}dt=i\int_{0}^{2\pi }\underset{ m=0}{\overset{\infty }{\sum
}}B_{mj}(L)(it)^{m}e^{-int}dt$

\begin{equation}
2\pi \cdot a_{nj}(L)=\int_{0}^{2\pi }\underset{m=0}{\overset{\infty }{\sum }} B_{mj}(L)(it)^{m}e^{-int}dt
\label{anj int}
\end{equation}

$\blacktriangle $ In order to show that we can exchange $\int_{0}^{2\pi }$and $\underset{\ \
m=0}{\overset{\infty }{\sum }}$in the above formula, it suffices to prove that $\int_{0}^{2\pi
}\underset{m=0}{\overset{\infty }{ \sum }}\left| B_{mj}(L)(it)^{m}e^{-int}\right| dt$ is finite.

$\int_{0}^{2\pi }\underset{m=0}{\overset{\infty }{\sum }}\left| B_{mj}(L)(it)^{m}e^{-int}\right|
dt=\int_{0}^{2\pi }\underset{\ \ m=0}{ \overset{\infty }{\sum }}\left| B_{mj}(L)\right| t^{m}dt\leq
\int_{0}^{2\pi } \underset{m=0}{\overset{\infty }{\sum }}\,\underset{k=-d}{\overset{d}{\sum }} \left|
a_{kj}(L)\right| \frac{k^{m}}{m!}t^{m}dt=\int_{0}^{2\pi }\underset{ k=-d}{\overset{d}{\sum }}\left|
a_{kj}(L)\right| \underset{\ \ m=0}{\overset{ \infty }{\sum }}\frac{(kt)^{m}}{m!}dt=\int_{0}^{2\pi }\underset{\
\ k=-d}{ \overset{d}{\sum }}\left| a_{kj}(L)\right| e^{kt}dt=\underset{\ \ k=-d}{ \overset{d}{\sum }}\left|
a_{kj}(L)\right| \int_{0}^{2\pi }e^{kt}dt$ which is finite.

$\blacktriangle $ Hence, we can exchange $\int_{0}^{2\pi }$and$\underset{m=0 }{\overset{\infty }{\sum }}$in
(\ref{anj int}), so $\ 2\pi \cdot a_{nj}(L)= \underset{m=0}{\overset{\infty }{\sum }}B_{mj}(L)\int_{0}^{2\pi
}(it)^{m}e^{-int}dt.$

Let $\lambda_{m,n}:=\frac{1}{2\pi }\int_{0}^{2\pi }(it)^{m}e^{-int}dt.$

\begin{equation}
\text{For any link }L,\text{ for any }n,j\in \mathbb{Z},\text{ }a_{nj}(L)= \underset{m=0}{\overset{\infty }{\sum
}}B_{mj}(L)\cdot \lambda _{m,n} \label{anj,sum}
\end{equation}

$\blacktriangle $ Let $v_{nj}^{N}(L)=\underset{m=0}{\overset{N}{\sum }} \lambda _{m,n}\cdot B_{mj}(L)$ and let
$\mu \in N^{\ast }.$ Let $n$, $j$ $ \in \mathbb{Z}$.

The restriction of $v_{nj}^{N}$ to links of $\mu $ components is $ v_{nj}^{N\mu
}(L)=\underset{m=0}{\overset{N}{\sum }}\lambda _{m,n}\cdot B_{mj}^{\mu }(L).$ It is a linear combination of the
Vassiliev invariants $ B_{mj}^{\mu }$ so it is a Vassiliev invariant. By (\ref{anj,sum}), for any
fixed $L\in \mathcal{L}^{(\mu )},\underset{N\rightarrow \infty }{\lim }%
v_{nj}^{N\mu }(L)=a_{nj}(L)$
\end{proof}

This proves that the restrictions of the coefficients of the HOMFLYPT
polynomial to links with a fixed number of component are pointwise limits of
Vassiliev invariants. By Proposition (\ref{WVCe}), they are pointwise limits
of Vassiliev invariants.

\begin{remark}
The approximation we get is totally explicit because we have the following formula to compute the coefficients
$\lambda _{m,n}=\underset{p=0}{\overset{ m-1}{\sum }}-\frac{\left( 2\pi i\right)
^{m-p-1}}{n}\frac{m!}{n^{p}.\left( m-p\right) !}$
\end{remark}

\section{Approximation of the coefficients of the Kauffman two-variable
polynomial by Vassiliev invariants}

We will use the Dubrovnik version of the Kauffman polynomial.

\begin{figure}[h]
\centerline{\psfig{figure=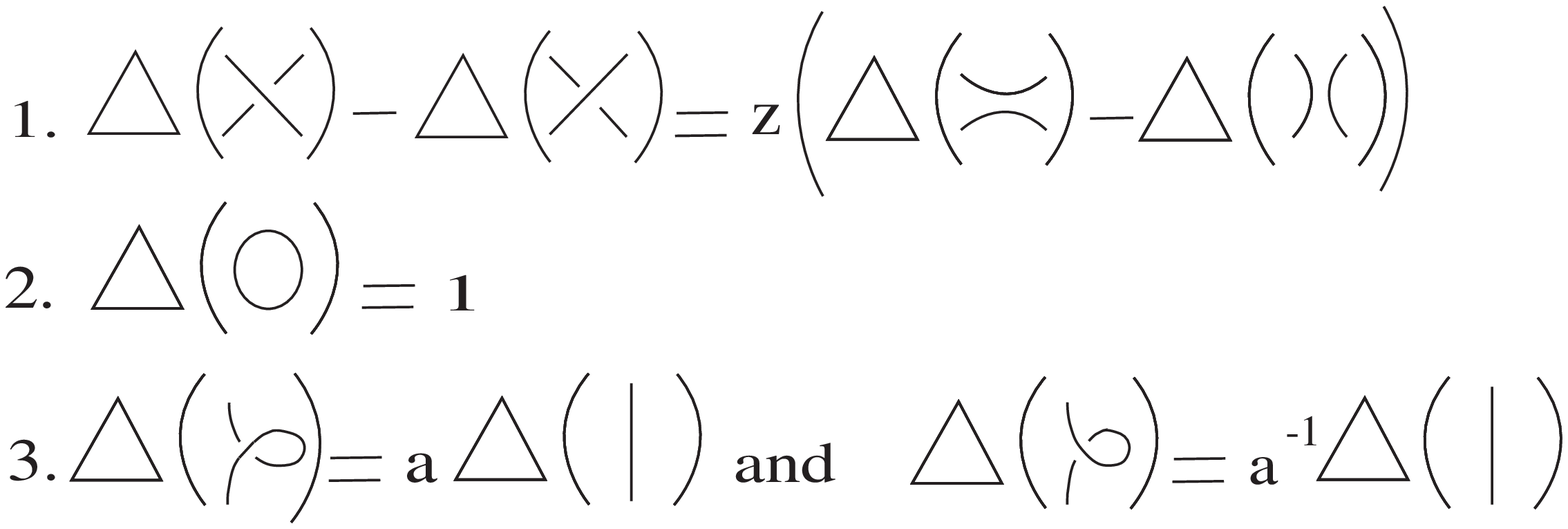,height=3.0cm}}
\caption{\textit{Axioms for $\Delta$}}
\label{DubrovnDltaDef}
\end{figure}

First, recall that $\Delta$ is defined by the axioms shown in Figure (\ref{DubrovnDltaDef}).

$\Delta $ is a regular isotopy invariant for unoriented links.\bigskip

Let $L$\ be an oriented link, let $\overrightarrow{D}$ be an oriented
diagram for $L$ and let $D$ be the corresponding unoriented diagram.

The \textit{Dubrovnik polynomial} is the knot invariant defined by $ F_{L}(a,z)=a^{-w(\overrightarrow{D})}\Delta
\left( D\right) $ where $w( \overrightarrow{D})$ is the Writhe number.\bigskip

The Dubrovnik polynomial is related to the Kauffman polynomial by the
following formula:

$F_{L}^{D}(a,z)=(-1)^{Comp(L)-1}F_{L}^{K}(a^{\prime }=ia,$ $z^{\prime }=-iz)$
where $F_{L}^{D}(a,z)$ is the Dubrovnik polynomial and $F_{L}^{K}(a^{\prime
},z^{\prime })$ is the Kauffman polynomial. Hence, once we prove that the
coefficients of the Dubrovnik polynomial are limits of Vassiliev invariants,
the same result will hold for the coefficients of the Kauffman
polynomial.\bigskip

The proof for the HOMFLYPT polynomial carries over to the Dubrovnik
polynomial with only slight modifications since we can get preliminary
results very similar to the ones we had for the HOMFLYPT polynomial, as
shown below.\bigskip

For instance, $F_{L}(a,z)\in \mathbb{Z}\left[ a^{\pm 1},z^{\pm 1}\right] $,
but as in the previous case, we can get a more precise statement.

\begin{lemma}
Let $\mu $ be the number of components of a link $L.$ Its Dubrovnik polynomial is a polynomial (not a Laurent
polynomial) in the variables $ a,a^{-1},$ $z$ and $\frac{a-a^{-1}+z}{z}$, the maximum possible power of $
\frac{a-a^{-1}+z}{z}$ being $\mu -1.$
\end{lemma}

\begin{proof}
The proof is similar to the proof of the corresponding result for the
HOMFLYPT polynomial.\medskip
\end{proof}

Also, we can use the same change of variable as in the case of the HOMFLYPT
polynomial to produce Vassiliev invariants:

\begin{lemma}
Let $W_{N,x}^{D}(L)$ be the polynomial obtained from the Dubrovnik polynomial of $L$ by setting $a:=e^{Nx}$ and
$z:=x.$ After power series expansion, each $w_{Nq}^{D}(L)$in $W_{N,x}^{D}(L)=\underset{q=0}{\overset{\infty
}{\sum }}w_{Nq}^{D}(L).x^{q}$ is a Vassiliev invariant of order (at most) $q$, for all $N$ $\in
\mathbb{Z}$.\medskip\
\end{lemma}

Since we used same change of variable as before, the power series expansion
will be the same and the proof is the same as the one for the HOMFLYPT
polynomial from now on.\medskip

\section{Acknowledgments}

I wish to thank Yongwu Rong for suggesting trying to prove this result, and
also for numerous helpful discussions, advice and support.

\end{document}